\documentclass[a4paper,oneside]{amsart}

\def\Wdn#1\wdn{\marginpar{\tiny #1}}
\long\def\WDN#1\wdn{[WDN: #1]\Wdn[Comment]\wdn}

\usepackage{amssymb}
\usepackage{amsmath}
\usepackage{amsfonts}
\usepackage{amsthm}
\usepackage{mathrsfs}
\usepackage{fix-cm}
\usepackage{graphicx}
\usepackage{amscd}

\usepackage{soul}
\sodef\spred{}{.2em}{.9em plus.4em}{1em plus.1em minus.1em}

\usepackage[latin1]{inputenc}
\usepackage[T1]{fontenc}

\usepackage[english]{babel}

\usepackage{bbm}

\usepackage[all]{xy}

\newbox\mybox
\def\overtag#1#2#3{\setbox\mybox\hbox{$#1$}\hbox to
  0pt{\vbox to 0pt{\vglue-#3\vglue-\ht\mybox\hbox to \wd\mybox
      {\hss$\ss#2$\hss}\vss}\hss}\box\mybox}
\def\undertag#1#2#3{\setbox\mybox\hbox{$#1$}\hbox to 0pt{\vbox to
    0pt{\vglue#3\vglue\ht\mybox\hbox to \wd\mybox
      {\hss$\ss#2$\hss}\vss}\hss}\box\mybox}
\def\lefttag#1#2#3{\hbox to 0pt{\vbox to 0pt{\vss\hbox to
      0pt{\hss$\ss#2$\hskip#3}\vss}}#1}
\def\righttag#1#2#3{\hbox to 0pt{\vbox to 0pt{\vss\hbox to
      0pt{\hskip#3$\ss#2$\hss}\vss}}#1}
\let\ss\scriptstyle

\def\Dot{\lower.2pc\hbox to 2.5pt{\hss$\bullet$\hss}}
\def\Circ{\lower.2pc\hbox to 2.5pt{\hss$\circ$\hss}}
\def\Vdots{\raise5pt\hbox{$\vdots$}}
\def\splicediag#1#2{\xymatrix@R=#1pt@C=#2pt@M=0pt@W=0pt@H=0pt}

\newcommand\lineto{\ar@{-}}
\newcommand\dashto{\ar@{--}}
\newcommand\dotto{\ar@{.}}

\newtheorem{thm}{Theorem}[section]

\newtheorem{prop}[thm]{Proposition}

\newtheorem{thm*}{Theorem}

\theoremstyle{definition} 
\newtheorem{defn}[thm]{Definition} 
\newtheorem{ex}[thm]{Example}

\newcommand{\Q}{\mathbbm{Q}}

\newcommand{\Z}{\mathbbm{Z}}

\newcommand{\num}[1]{\lvert #1 \rvert}

\newcommand{\sign}{\operatorname{sign}}
\newcommand{\lcm}{\operatorname{lcm}}

\DeclareMathOperator*{\conect}{\#}

\newcommand{\inv}{^{-1}}

\makeatletter
\newcommand\@b@gconect[1]{%
\vcenter{\hbox{#1$\m@th\mkern2mu\conect\mkern2mu$}}}
\newcommand\@bigconect{%
\mathchoice{\@b@gconect\huge} 
{\@b@gconect\LARGE} 
{\@b@gconect{}} 
{\@b@gconect\footnotesize} 
}
\newcommand\bigconect{\mathop{\@bigconect}\displaylimits}
\makeatother

\renewcommand{\phi}{\varphi}
\renewcommand{\epsilon}{\varepsilon}

\begin{document}

\bibliographystyle{alpha}


\title[Construction universal abelian covers]
{Constructing universal abelian covers of graph manifolds}
\author{Helge M\o{}ller Pedersen}
\address{Matematische Institut\\ Universität Heidelberg
\\ Heidelberg, 69120}
\email{pedersen@mathi.uni-heidelberg.de}
\keywords{rational homology sphere,
abelian cover}
\subjclass[2000]{57M10, 57M27}
\begin{abstract}
  To a rational homology sphere graph manifold one can associate a
  weighted tree invariant called splice diagram. It was shown in
  \cite{myarticle} that the splice diagram determines the universal
  abelian cover of the manifold. We will in this article turn the
  proof of this in to an algorithm to explicitly construct the
  universal abelian cover from the splice diagram.
\end{abstract}
\maketitle

\section{Introduction}

Graph manifolds is an important class of $3$-manifolds, they are
defined as the manifolds which have only Seifert fibered pieces in their
JSJ-decomposition. They are also the $3$-manifolds which are
boundaries of plumbed $4$-manifolds, and a very used method to
represent a graph manifold $M$ is by giving a plumbing diagram of a
$4$-manifold $X$ such that $M=\partial X$. Neumann gave a complete
calculus for changing $X$ but keeping $M$ fixed in \cite{plumbing},
and when we are going to construct the universal abelian cover in
section \ref{construction} we are going to do this by constructing a plumbing
diagram of it. It should also be noted that graph manifolds are also
the manifolds that have no hyperbolic pieces in their geometric
decomposition.

If one restrict to rational homology spheres (from know
on $\Q$HS) then one have a graph invariant of graph manifolds called
splice diagrams. They where original introduced in \cite{EisenbudNeumann} and
\cite{Siebenmann} for integer homology sphere graph manifolds, and were
then later generalized by Neumann and Wahl to $\Q$HS's in
\cite{NeumannWahl3}, and used extensively in
\cite{neumannandwahl1} and \cite{neumannandwahl2}.

The simplest example of a graph manifold is of course a Seifert fibered
manifold, and if one restrict to $\Q$HS's, Neumann found a nice
construction of the universal abelian 
cover in \cite{neumann83a}, namely as the link of a Brieskorn complete
intersection defined 
by the collection of the first of the two Seifert invariants
associated to the singular fibers. This is exactly the information
given by the splice diagram of Seifert fibered manifolds, and Neumann
and Wahl used the splice diagram to generalize the Brieskorn complete
intersections provided the splice diagram satisfy what they called the
semigroup condition, to what they called \emph{splice diagram
  equations}. Under a further restriction on the given manifold they
were able to prove that the link of the splice diagram equations is the
universal abelian cover in \cite{neumannandwahl2}. 

This indicates that the splice diagram might determine the universal
abelian cover, and in \cite{myarticle} I was able to prove the
following theorem:
\begin{thm}\label{univeralsabcover}  
Let $M_1$ and $M_2$ be two $\Q$HS graph manifold there have the same
splice diagram. Let $\widetilde{M}_i\to M_i$ be the universal abelian
covers. Then $\widetilde{M}_1$ and $\widetilde{M}_2$ are homeomorphic.
\end{thm}
The proof consist of inductively constructing the universal abelian
cover from the splice diagram, and the purpose of this article is to
extract an algorithm for constructing the universal abelian cover from
the proof. I will hence not prove that I actually construct the
universal abelian cover, but refer to the proof given in
\cite{myarticle} for Theorem \ref{univeralsabcover}, (in that article
it is Theorem 6.3).

Returning to splice diagram singularities, i.e.\ the complete
intersections define by the splice diagram equations of a splice
diagram $\Gamma$, then one can use
this algorithm to construct a dual resolution diagram, provided that
there is a manifold (orbifold) there satisfies the (orbifold)
congruence condition and has $\Gamma$ as its splice diagram. This is
for example always true if $\Gamma$ has only two nodes see
\cite{orbifoldcongruencecondition}.   

The algorithm is as mention above going to be give in section
\ref{construction}, section \ref{splicediagram} will introduce splice
diagrams, their relation with plumbing diagrams and mention some
results needed for the algorithm. 

\section{Splice Diagrams}\label{splicediagram} 

A \emph{splice diagram} is a weighted tree with no vertices of valence
two. At vertices of valence greater than two, we call such vertices
for \emph{nodes}, one adds a sign, and on edges adjacent to nodes on
adds a non negative integer weights. 

To any $\Q$HS graph manifold $M$ we can associate a splice diagram
$\Gamma(M)$ by the following procedure:
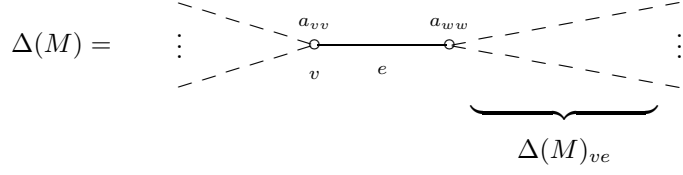
\begin{itemize}
\item Take a node (a vertex which are going to end up as a node) for
  each Seifert fibered piece of the JSJ-decomposition of $M$, we will
  not distinguish between the nodes and the corresponding Seifert
  fibered pieces.
\item Connect two nodes if they are glued in the JSJ-decomposition to
  create $M$.
\item Add a \emph{leaf} (a valence one vertex connected by an edge) to a node
  for each Singular fiber of the Seifert fibration.
\item Adds the sign of the linking number of two nonsingular fibers at
  a node. See \cite{myarticle} for how to define linking numbers.
\item Let $v$ be a node and $e$ an edge adjacent to $v$. Then the edge
  weight $d_{ve}$ is determined the following way. Cut $M$ along the
  torus $T$ corresponding to $e$ (either a torus from the
  JSJ-decomposition of $M$ or the boundary of a tubular neighborhood of
  a singular fiber) into the pieces $M_v$ and $M_{ve}'$, where $v$ is
  in $M_v$. Then glue a solid torus into the boundary of $M_{ve}'$ by
  identifying a meridian with the image of a fiber from $M_v$, and
  call this new closed graph manifold $M_{ve}$. Then
  $d_{ve}=\num{H_1(M_{ve})}$ if $H_1(M_{ve})$ is finite or $0$ otherwise.
\end{itemize}

The standard way to represent graph manifolds are by plumbing
diagrams, and we will next describe how to get the splice diagram from
a plumbing diagram $\Delta$ of $M$.

To construct the graph structure of $\Gamma(M)$ from $\Delta$ on just
suppress all vertices of valence two, i.e.\ replacing any configuration
like   
$$
\xymatrix@R=6pt@C=24pt@M=0pt@W=0pt@H=0pt{
\undertag{\Circ}{v}{4pt}\lineto[r] & \overtag{\Circ}{-b_1}{8pt}\lineto[r]&
\overtag{\Circ}{-b_2}{8pt}\dashto[rr]&&
\overtag{\Circ}{-b_k}{8pt}\lineto[r] &\undertag{\Circ}{w}{4pt} 
\hbox to 0pt {\hss}}
$$
with an edge
$$
\xymatrix@R=6pt@C=24pt@M=0pt@W=0pt@H=0pt{
\undertag{\Circ}{v}{4pt}\lineto[r]  &\undertag{\Circ}{w}{4pt}
\hbox to 0pt {~.\hss}}
$$
Let $A(\Delta)$ be the intersection matrix of the $4$ manifold defined
by $\Delta$.

The edge weights and signs are found by the following propositions from
\cite{myarticle}.
\begin{prop}
Let $v$ be a node in $\Gamma(M)$, and $e$ be an edge on that node. We
get the weight $d_{ve}$ on that edge by
$d_{ve}=\num{\det(-A(\Delta(M)_{ve}))}$, where $\Delta(M)_{ve}$ is is the
connected component of $\Delta(M)-e$ which does not contain $v$.
\end{prop}

$$
\xymatrix@R=6pt@C=24pt@M=0pt@W=0pt@H=0pt{
&&&&&&&\\
\Delta(M)=&\Vdots&&\undertag{\overtag{\Circ}{a_{vv}}{8pt}}{v}{8pt}
\lineto[rr]\dashto[ull]\dashto[dll]&{}\undertag{}{e}{8pt}&
\overtag{\Circ}{a_{ww}}{8pt}\dashto[urr]\dashto[drr]&&\Vdots\\
&&&&&&&\\
&&&&&&{\hbox to 0pt{\hss$\underbrace{\hbox to 70pt{}}$\hss}}&\\
&&&&&&{\Delta(M)_{ve}}&}$$


\begin{prop}
Let $v$ be a node in $\Gamma(M)$. Then the sign $\epsilon$ at $v$ is
$\epsilon=-\sign(a_{vv})$, where $a_{vv}$ is the entry of 
$A(M)\inv$ corresponding to the node $v$. 
\end{prop}

The most important information the splice diagram do not encode of the
manifold is $\num{H_1(M)}$, but that and the splice diagram do
determine the rational euler number of any of the Seifert fibered
pieces of $M$, by the following proposition from \cite{myarticle}.

\begin{prop}\label{eulernumber}
Let $v$ be a node in a splice diagram decorated as in Fig.\ 1 below
with $r_i\neq 0$ for $i\neq 1$, and let $e_v$ be the rational euler
number of $M_v$. Then
\begin{align}  
e_v=-d\big(\frac{\epsilon s_1}{ND_1\prod_{j=2}^kr_k}+
\sum_{i=2}^k\frac{\epsilon_iM_i}{r_iD_i}\big) 
\end{align}
where $d=\num{H_1(M)}$, $N=\prod_{j=1}^kn_j$,
$M_i=\prod_{j=1}^{l_i}m_{ij}$, and $D_i$ is the edge determinant associated to
the edge between $v$ and $v_i$.
\end{prop}
$$\splicediag{8}{30}{
  &&&&\\
 &&& \overtag\Circ{v_1}{8pt}
      \lineto[ur]^(.5){m_{11}}
  \lineto[dr]_(.5){m_{1l_l}}&\Vdots\\
  \Circ&&&&\\ 
  \Vdots&\overtag\Circ {v} {8pt}\lineto[ul]_(.5){n_1}
  \lineto[dl]^(.5){n_k}
  \lineto[uurr]^(.25){r_1}^(.75){s_1}
  \lineto[ddrr]^(.25){r_k}^(.75){s_k}
  &\Vdots&&  \\
  \Circ&&&&\\
  &&& \overtag\Circ{v_k}{8pt}
      \lineto[ur]^(.5){m_{k1}}
  \lineto[dr]_(.5){m_{kl_k}}&\Vdots\\
  &&&&\\
  &&\text{Figure 1}&&\hbox to 0 pt{~.\hss} }$$
Note that this does give a formula for $e_v/d$ from $\Gamma$, which we
will need later.

In the algorithm for constructing the universal abelian cover of $M$
form $\Gamma(M)$, one number associated to each end of an edge in
$\Gamma(M)$ is going to be very important, the \emph{ideal
  generator}, which is constructed the following way. Let $v$ and $w$
be two vertices of $\Gamma(M)$ then we define the \emph{linking
  number} of $v$ and $w$ $l_{vw}$ as the product of all edge weights
adjacent to but not on the shortest path from $v$ to $w$. We define
$l'_{vw}$ the same way, but omitting weights adjacent to $v$ and $w$.
If $e$ is an edge adjacent to $v$, we then let $\Gamma_{ve}$
be the connected component of $\Gamma(M)-e$ not containing $v$. And
define the following ideal of $\Z$
\begin{align*}
I_{ve}&=\langle l'_{vw}\vert\ w \text{ a leaf in } \Gamma_{ve}\rangle.
\end{align*} 
Then we define the ideal generator $\overline{d}_{ve}$ associated to
$v$ and $e$ to be the positive generator of $I_{ve}$. 

\begin{defn}
A splice diagram $\Gamma$ satisfy the \emph{ideal condition} if the
ideal generator $\overline{d}_{ve}$ divides the edge weight $d_{ve}$
\end{defn}

\begin{prop}
Let $M$ be a $\Q$HS graph manifold, then $\Gamma(M)$ satisfy the ideal
condition.
\end{prop}

This proposition follows from the following topological description of
the ideal generator from Appendix 1 of \cite{neumannandwahl2}.

\begin{thm}
The ideal generator $\overline{d}_{ve}$ is $\num{H_1(M_{ve}/K)}$, where
  $K$ is the knot given as the core of the solid torus glued into
  $M_{ve}'$ to construct $M_{ve}$.
\end{thm}

\section{Construction the Universal Abelian Cover: An Example}\label{construction}

In this section we are going to see how the algorithm used in the
proof of Theorem \cite{myarticle} can be used to construct the
universal abelian cover $\widetilde{M}$ of a graph manifold $M$, from the splice
diagram $\Gamma(M)$. We are going to specify $\widetilde{M}$ by
constructing a plumbing diagram $\Delta$ for $\widetilde{M}$. To
illustrate the construction we are going to use the following example

$$\splicediag{8}{30}{
  &&&&&\\
  &\Circ &&&&\Circ \\
  \Gamma=&&\overtag\Circ {v_0} {8pt}\lineto[ul]_(.5){3}
  \lineto[dl]^(.5){18}
  \lineto[rr]^(.25){23}^(.75){15}&& \overtag\Circ{v_1}{8pt}
  \lineto[ur]^(.5){2}
  \lineto[dr]_(.5){3}&\\
  &\Circ &&&&\Circ \\
&&&&&\hbox to 0 pt{~.\hss} }$$
There are four different manifolds who has $\Gamma$ as their splice
diagram, and several more non manifold graph orbifolds. By Theorem
\cite{myarticle} $\Gamma$ is the splice diagram of a singularity
link, it then follows from \cite{rationalcover} that $\widetilde{M}$ is
a rational homology sphere. The example is also interesting since non
of the manifolds with $\Gamma$ as their splice diagram satisfy the
congruence condition of Neumann and Wahl see \cite{neumannandwahl2},
but there are non manifold orbifolds with splice diagram $\Gamma$ which
satisfy the orbifold congruence condition see
\cite{orbifoldcongruencecondition}. Below is plumbing diagrams for the
four manifolds with $\Gamma$ as their splice diagram:

$$
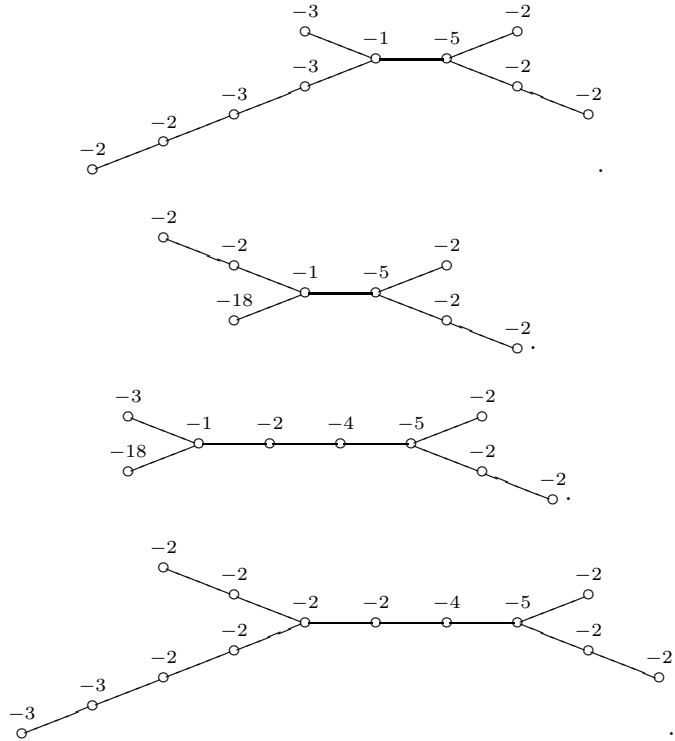
\xymatrix@R=6pt@C=24pt@M=0pt@W=0pt@H=0pt{
&&&\overtag{\Circ}{-3}{8pt}  &&&\overtag{\Circ}{-2}{8pt}&\\
&&&& \overtag{\Circ}{-1}{8pt}\lineto[ul]\lineto[dl]\lineto[r] &
\overtag{\Circ}{-5}{8pt}\lineto[ur]\lineto[dr]&&\\
&&&\overtag{\Circ}{-3}{8pt}\lineto[dl]&&&
\overtag{\Circ}{-2}{8pt}\lineto[dr]&\\ 
&&\overtag{\Circ}{-3}{8pt}\lineto[dl]&&&&&\overtag{\Circ}{-2}{8pt}\\
&\overtag{\Circ}{-2}{8pt}\lineto[dl]&&&&&&\\
\overtag{\Circ}{-2}{8pt} &&&&&&&\hbox to 0pt {~.\hss}}
$$

$$
\xymatrix@R=6pt@C=24pt@M=0pt@W=0pt@H=0pt{
\overtag{\Circ}{-2}{8pt} &&&&&\\
&\overtag{\Circ}{-2}{8pt}\lineto[ul]  &&&\overtag{\Circ}{-2}{8pt}&\\
&& \overtag{\Circ}{-1}{8pt}\lineto[ul]\lineto[dl]\lineto[r] &
\overtag{\Circ}{-5}{8pt}\lineto[ur]\lineto[dr]&&\\
&\overtag{\Circ}{-18}{8pt}&&&
\overtag{\Circ}{-2}{8pt}\lineto[dr]&\\ 
&&&&&\overtag{\Circ}{-2}{8pt} \hbox to 0pt {~.\hss}}
$$

$$
\xymatrix@R=6pt@C=24pt@M=0pt@W=0pt@H=0pt{
\overtag{\Circ}{-3}{8pt}  &&&&&\overtag{\Circ}{-2}{8pt}&\\
& \overtag{\Circ}{-1}{8pt}\lineto[ul]\lineto[dl]\lineto[r] &
\overtag{\Circ}{-2}{8pt}\lineto[r]& \overtag{\Circ}{-4}{8pt}\lineto[r] &
\overtag{\Circ}{-5}{8pt}\lineto[ur]\lineto[dr]&&\\
\overtag{\Circ}{-18}{8pt}&&&&&
\overtag{\Circ}{-2}{8pt}\lineto[dr]&\\ 
&&&&&&\overtag{\Circ}{-2}{8pt} \hbox to 0pt {~.\hss}}
$$

$$
\xymatrix@R=6pt@C=24pt@M=0pt@W=0pt@H=0pt{
&&\overtag{\Circ}{-2}{8pt} &&&&&&&\\
&&&\overtag{\Circ}{-2}{8pt}\lineto[ul] &&&&&\overtag{\Circ}{-2}{8pt}&\\
&&&& \overtag{\Circ}{-2}{8pt}\lineto[ul]\lineto[dl]\lineto[r] &
\overtag{\Circ}{-2}{8pt}\lineto[r]& \overtag{\Circ}{-4}{8pt}\lineto[r] &
\overtag{\Circ}{-5}{8pt}\lineto[ur]\lineto[dr]&&\\
&&&\overtag{\Circ}{-2}{8pt}\lineto[dl]&&&&&
\overtag{\Circ}{-2}{8pt}\lineto[dr]&\\ 
&&\overtag{\Circ}{-2}{8pt}\lineto[dl] &&&&&&&\overtag{\Circ}{-2}{8pt}\\
&\overtag{\Circ}{-3}{8pt}\lineto[dl] &&&&&&&&\\
\overtag{\Circ}{-3}{8pt} &&&&&&&&&\hbox to 0pt {~.\hss}}
$$

\subsection{Constructing the Building Blocks}

The inductive procedure in the construction of the universal abelian
cover works by taking an edge $e$ between two nodes of $\Gamma$ and make
a new non connected splice diagram $\Gamma_e$ where $e$, has been
replaced with two leaves. So starting with the edge called $e_1$ and
going going through this process of cutting the edge until we have cut
the last edge between two nodes $e_{N-1}$, we get that $\Gamma_{e_{N-1}}$ is a
collection of one node splice diagrams
$\Gamma_{e_{N-1}}=\{\Gamma_i\}_{i=1}^N$. For 
  each of thees one node splice diagram $\Gamma_i$ one then takes a
  number of copies of a specific manifold $M_i$, and use the
  information from  the $\Gamma_{e_j}$'s to glue the pieces
together. So the first step is to determine this manifolds
$\{M_i\}_{i=1}^N$, which are the building blocks of the universal
abelian cover.  

First lets see how the $\Gamma_{e_j}$'s are going to
look. Each time we cut an edge $e$ between the nodes $w_1$ and $w_2$ in
$\Gamma$, we divide every edge 
weight $d_{ve'}$ such that $w_1$ or $w_2$ is in $\Gamma_{ve'}$, by the
ideal generator 
$\overline{d}_{w_ie}$ of the edge weight $d_{w_ie}$ such that $v$ is not
in $\Gamma_{w_ie}$.
 In our example we only have two edge weights where this is true with
 respect to $e$  namely
$d_{v_0e}=23$ and $d_{v_1e}=15$, and $\overline{d}_{v_0e}=1$ and
$\overline{d}_{v_1e}=3$. So the two one node splice diagrams
$\Gamma_1$ and $\Gamma_2$ are going to look like 

$$\splicediag{8}{30}{
  &  \Circ &&&&&&\Circ \\
  \Gamma_1=&&\overtag\Circ {v_0} {8pt}\lineto[ul]_(.5){3}
  \lineto[dl]^(.5){18}
  \lineto[r]^(.25){23}& \overtag\Circ {(1,1)} {8pt},&\Gamma_2=&  \overtag\Circ
  {(1,3)} {8pt}\lineto[r]^(.75){5}& \overtag\Circ{v_1}{8pt}
  \lineto[ur]^(.5){2}
  \lineto[dr]_(.5){3}&\\
  &\Circ &&&&&&\Circ\hbox to 0 pt{~.\hss} }$$
The pair added to the new leaves is recording of the following
information which is going to be used when the gluing are made: the
first number specifies which number in the sequence of cutting this
is, in this case the first, and the second number is the
ideal generator associated to the weight before cutting.     

Next we want to find the building block $M_i$ associated to each of
the $\Gamma_i$'s. To do this we have to separate the $\Gamma_i$'s into
two types, the first are the once that do not have an edge weight of
$0$, and the second are the once that, remember at most one
weight adjacent to a node can be $0$.

 In the first case we use the
following theorem
\begin{thm}
Let $M$ be a rational homology orbifold fibration $S^1$-fibration over
a orbifold surface, with Seifert invariants
$(\alpha_1,\beta_1),\dots,(\alpha_n,\beta_n)$. Then the universal
abelian cover of $M$ is the link of the Brieskorn complete intersection
$\Sigma(\alpha_1,\dots,\alpha_n)$.
\end{thm}
The way one construct the manifolds after cutting an edge may result
in graph orbifolds instead of just graph manifolds, as explained in
the proof of 6.3 in \cite{myarticle}, and hence we need this theorem
for $S^1$ orbifold fibrations.
Neumann only proves this theorem for Seifert fibered manifolds in
\cite{neumann83a} and \cite{neumann83b}, but the proof given in
\cite{neumann83b} also work in the general case of an orbifold
$S^1$-fibration. The value of $\epsilon$ does not matter, since
reversing the orientation of a Seifert fibered manifold only changes
the $\beta_i$'s not the $\alpha_i$'s and hence only change the splice
diagrams by replacing $\epsilon$ with $-\epsilon$.  

So in our example $M_1$ is the link of $\Sigma(3,18,23)$
and $M_2$ is the link of $\Sigma(2,3,5)$. 

Next we use the following theorem to get plumbing diagrams for the
$M_i$'s.
\begin{thm}
Let $M$ be the link of the Brieskorn complete intersection
$\Sigma(\alpha_1,\dots,\alpha_n)$. A plumbing diagram for $M$ is
given by
$$
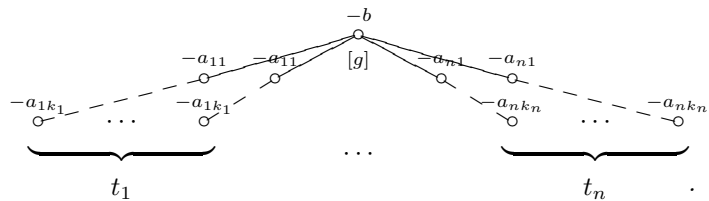
\xymatrix@R=6pt@C=24pt@M=0pt@W=0pt@H=0pt{
&&&&\undertag{\overtag{\Circ}{-b}{8pt}}{[g]}{4pt} \lineto[ddll]
\lineto[ddl] \lineto[ddr] \lineto[ddrr]&&&&\\
&&&&&&&&\\
&&\overtag{\Circ}{-a_{11}}{8pt}\dashto[ddll]
&\overtag{\Circ}{-a_{11}}{8pt}\dashto[ddl] &  &
\overtag{\Circ}{-a_{n1}}{8pt}\dashto[ddr] &
\overtag{\Circ}{-a_{n1}}{8pt}\dashto[ddrr] &&\\
&&&&&&&&\\
\overtag{\Circ}{-a_{1k_1}}{8pt} &\cdots
&\overtag{\Circ}{-a_{1k_1}}{8pt} & &&& \overtag{\Circ}{-a_{nk_n}}{8pt}
& \cdots & \overtag{\Circ}{-a_{nk_n}}{8pt}\\
&{\hbox to 0pt{\hss$\underbrace{\hbox to 70pt{}}$\hss}}
&&&\cdots&&&{\hbox to 0pt{\hss$\underbrace{\hbox to 70pt{}}$\hss}}& \\
& {t_1} &&&&&& {t_n}& \hbox to 0pt {~.\hss}}
$$
The values of $g$, and the $t_i$'s, are given by 
\begin{align}
t_i &= \frac{\prod_{j\neq i}(\alpha_j)}{\lcm_{j\neq i}(\alpha_j)} \\
g
&=\tfrac{1}{2}\big(2+\frac{(n-2)\prod_i\alpha_i}{\lcm_i(\alpha_i)}-\sum_{i=1}^n t_i\big).
\end{align}
Then one calculates numbers $p_1,\dots,p_n$ as 
\begin{align}
p_i &= \frac{\lcm_j(\alpha_j)}{\lcm_{j\neq i}(\alpha_j)},
\end{align}
and find numbers $q_1,\dots, q_n$ as the smallest possible solutions
to the equations
\begin{align}
\frac{\lcm_j(\alpha_j)}{\alpha_i}q_i &\equiv -1 (\mod p_i).
\end{align}
Then the $a_{ij}$'s, are given as the continued fraction
$p_i/q_i=[a_{i1},\dots,a_{ik_i}]$. Finally $b$ is given by
\begin{align}   
b&= \frac{\prod_i\alpha_i+\lcm_i(\alpha_i)\sum_iq_i\prod_{j\neq
    i}\alpha_j}{(\lcm_i\alpha_i)^2}.
\end{align}
\end{thm}

Before we use this theorem to make a plumbing diagram $\Delta_i$ for the $M_i$,
notice that to make the gluing we have to remove some solid tori from
the $M_i$'s to make the gluing, so we need to record this data in $\Delta_i$.
 Some leaves in $\Gamma_i$ have a pair of integers
attached. These leaves corresponds to the tori in $M$ we cut along
when we created $\Gamma_i$. Since $M_i$ is the universal abelian
cover of any graph orbifold with $\Gamma_i$ as its splice diagram,
several singular fibers sits above the singular fiber corresponding to
these leaves. It is a neighborhood of each of these singular
fibers we have to remove. So if $\alpha_j$ is a edge weight in
$\Gamma_i$ to a leaf with a pair attached, the the $t_j$ singular
fibers above the leaf, corresponds to all the strings with the weights
$-a_{j1},\dots,-a_{jn_j}$. So in the plumbing diagram for $M_i$ we
replace these string with an arrow, and add a triple which consist of 
the pair attached to the leaf and $p_j/q_j=[a_{j1},\dots,a_{jk_j}]$ to
each of the arrows. 

Using these theorems on our example we get the following plumbing
diagrams
$$
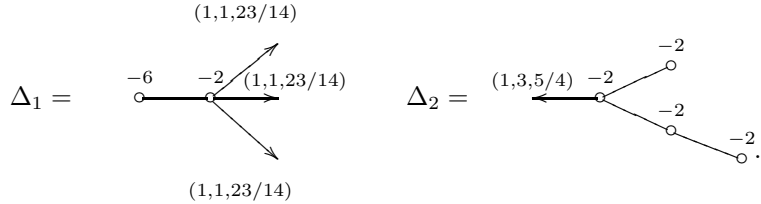
\xymatrix@R=6pt@C=24pt@M=0pt@W=0pt@H=0pt{
&&&&&&&&&\\
&&&&&&&&\overtag{\Circ}{-2}{8pt}&\\
\Delta_1=& \overtag{\Circ}{-6}{8pt}\lineto[r]
&\overtag{\Circ}{-2}{8pt}\ar[uur]^(1.25){(1,1,23/14)}
\ar[r]^(1.25){(1,1,23/14)}
\ar[ddr]_(1.25){(1,1,23/14)} &&&\Delta_2= && \overtag{\Circ}{-2}{8pt}
\lineto[ur]\ar[l]_(1){(1,3,5/4)}
\lineto[dr] &&\\
&&&&&&&&\overtag{\Circ}{-2}{8pt} \lineto[dr] &\\
&&&&&&&&&\overtag{\Circ}{-2}{8pt}
\hbox to 0pt {~.\hss}}
$$  

The second case is not as easy, the proof of Theorem 6.3 in
\cite{myarticle} give a construction in this case, but it might not be
a Seifert fibered manifold, an I have at the present no simple way to
find a plumbing diagram for the building blocks in this case.


\subsection{Gluing the Building Blocks}

The only thing that remains to construct the universal abelian cover
is to glue together the building blocks $M_i$, this will be done by
using the plumbing diagrams $\Delta_i$ to create a plumbing diagram
$\Delta$ for $\widetilde{M}$.  

Start by taking two of the $\Delta_i$'s and create
a plumbing diagram $G_1$, then we take an other of the $\Delta_i$'s
and glue this to $G_1$ to create $G_2$, we continue this process until
all the $\Delta_i$'s has been used, and then $\Delta=G_{N-1}$ where
$G_{N-1}$ is the last created plumbing diagram. 

Now the order we glue the $\Delta_i$'s together in is important, this
is why we added a triple at the arrows. We start by taking $\Delta_i$
which have at least one arrow that has a triple $(N-1,d_i,r_i)$, where
$N-1$ is the highest value for value for the first number in the
triple. And that $\Delta_j$ such at least one arrow that has a triple
$(N-1,d_j,r_j)$. By the method we constructed the $\Delta_i$'s
there are exactly two satisfying this. Then we take $d_i$ copies of
$\Delta_i$ and $d_j$ copies of $\Delta_j$. We the create an
intermediate $\widetilde{G}_1$ by for each of the copies of $\Delta_i$
replace each arrow with the triple $(N-1,d_i,r_i)$ with a dashed line
to a copy of $\Delta_j$ replacing the arrow with the triple
$(N-1,d_j,r_j)$, such that a copy of $\Delta_i$ is only connected to a
copy of $\Delta_j$ once. This will create a connected weighted graph
$\widetilde{G}_1$, with no arrows which has first number in the triple
equal to $N-1$.

Lets see how this is done in our example. We only have
two $\Delta_i$'s, so we start by gluing $\Delta_1$ to $\Delta_2$. The
triples are $(1,1,23/14)$ and $(1,3,5/4)$ so we start by taking one
copy $\Delta_1$ and 3 copies of $\Delta_2$, replace each of the arrows
in the copy of $\Delta_1$ with a dashed line to one of the copies of
$\Delta_2$ replacing its arrow. We then get 
$$
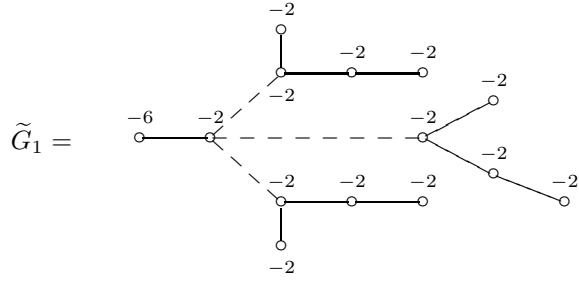
\xymatrix@R=6pt@C=24pt@M=0pt@W=0pt@H=0pt{
&&&&&&&\\
&&&\overtag{\Circ}{-2}{8pt}&&&&\\
&&&&&&&\\
&&&\undertag{\Circ}{-2}{4pt}\lineto[uu]\lineto[r] 
& \overtag{\Circ}{-2}{8pt}\lineto[r]&\overtag{\Circ}{-2}{8pt}
&&\\
&&&&&&
\overtag{\Circ}{-2}{8pt}&\\ \widetilde{G}_1=&
\overtag{\Circ}{-6}{8pt}\lineto[r] 
&\overtag{\Circ}{-2}{8pt}\dashto[uur] \dashto[rrr]
\dashto[ddr] &&& \overtag{\Circ}{-2}{8pt}
\lineto[ur]
\lineto[dr] &&\\
&&&&&& \overtag{\Circ}{-2}{8pt} \lineto[dr] &\\
&&&\overtag{\Circ}{-2}{8pt}\lineto[dd]\lineto[r] & 
\overtag{\Circ}{-2}{8pt}\lineto[r] & \overtag{\Circ}{-2}{8pt}
&& \overtag{\Circ}{-2}{8pt}\\
&&&&&&&\\
&&&\undertag{\Circ}{-2}{4pt}&&&&\\
&&&&&&&\\
&&&&&&&\hbox to 0pt {~.\hss}}
$$ 

The next step is to replace the dashed lines by a string to create
$G_1$. First by symmetry all the strings are going to be equal, so we
only have to calculate one string. To do this we use that there are
two different ways to calculate the rational euler number of the
Seifert fibered piece corresponding to a node in $G_1$, one using
$G_1$ and one given by the splice diagram by a formula derived at the
end of the proof of Theorem 6.3 in \cite{myarticle}. 

Choose a node $v$ of $\widetilde{G}_1$ which is attached to a dashed
line, then the rational euler number is given by $b+\sum_{e} q_e/p_e$,
where the sum is taking over all edges adjacent to $v$ (including the
dashed lines), $(p_e,q_e)$ is the Seifert pair associated to the
string and $b$ is the weight at $v$. Now there are four types of
different edges attached to $v$ 
and we need to see how to get $(p_e,q_e)$ from each type of the edge. 

First there are the edge that starts a string that ends at a valence
one vertex. Form these strings we get get $(p_e,q_e)$ from the
continued fraction associated to the string i.e.\
$p_e/q_e=[a_{e1},\dots,a_{ek_e}]$.  

Second types of edges are on string that leads to other nodes (when
one makes $G_1$ these do not exist, but they can be there when we are
going to make $G_2$). We again gain the Seifert pair from the
continued fraction, this time from the string between $v$ and the
other node.

The third type of edges are the arrows, their we gain $(p_e,q_e)$ from
the triple $(n_e,d_e,r_e)$ attached to the arrow as $p_e/q_e=r_e$. 

The last type of edges are the dashed lines, here we do not find the
Seifert pair since we are trying to make an equation to do just
that. But notice that all the dashed lines have the same Seifert
pair $(p,q)$, hence we get the following equation
\begin{align}\label{eq1}
e_v &=d'\frac{q}{p}+b+\sum_e\frac{q_e}{p_e}.
\end{align}
Where the sum is taken over all edges at $v$ except the dashed lines,
and $d'$ is the number of dashed lines at $v$. Notice that if $v$ is a
node sitting above the $\Delta_i$ piece, then $d'=d_j$. 

Returning to our example, if we use the left most node as $v$, the
equation becomes
\begin{align*}
e_v&=3\frac{p}{q}-2+\frac{1}{6}=3\frac{p}{q}-\frac{11}{6}.
\end{align*}

Now $G_1\bigcup(\bigcup_{l\neq i,j}\Delta_l)$ is a plumbing diagram for
any manifold with splice diagram $\Gamma_{e_{N-2}}$, hence it is
$\Gamma_{e_{N-2}}$ we need to use when we make the other calculation
of $e_v$.

From the end of the proof of Theorem 6.3 in \cite{myarticle} one gets
that if $v$ sits above the $\Delta_i$ piece
\begin{align}\label{eq2}
e_v &=\frac{\lambda^2}{d_i}\tilde{e}_v/d.
\end{align}
Where $\lambda=\prod
m_j/\lcm(m_1/\overline{d}_1,\dots,m_k/\overline{d}_k)$ where the $m_j$'s
are the edge weight adjacent to the node corresponding to $v$ in
$\Gamma_{e_{N-2}}$ and the $\overline{d}_j$'s are the ideal generator
associated to the edges. $\tilde{e}_v$ is the rational euler number of
the the node corresponding to $v$ in any graph orbifold
$\widetilde{M}$ whit $\Gamma(\widetilde{M})=\Gamma_{e_{N-2}}$ and
$d=\num{H_1^{orb}(\widetilde{M})}$, now neither of these numbers are
determined by $\Gamma_{e_{N-2}}$, but proposition \ref{eulernumber}
gives a formula for $\tilde{e}_v/d$ only using $\Gamma_{e_{N-2}}$. 

In our example $\Gamma_{e_{N-2}}=\Gamma$, so using this we find that
$\lambda=3$ and $\tilde{e}_v/d=-5/378$, so $e_v=-5/42$.

Now one find $p/q$ by combining the equations \eqref{eq1} and
\eqref{eq2}, which in our example gives $p/q=4/7$. Remember this is
the continued fraction associated to the string replacing the dashed
lines when seen from $v$, if we have used the node in the other end of
the line we would have found $p/q'$ where $qq'\equiv -1(\mod
p)$. Replacing all the dashed lines with the strings corresponding to
the continued fractions on gets $G_1$ from $\widetilde{G}_1$. This
becomes the following plumbing diagram in our example
$$
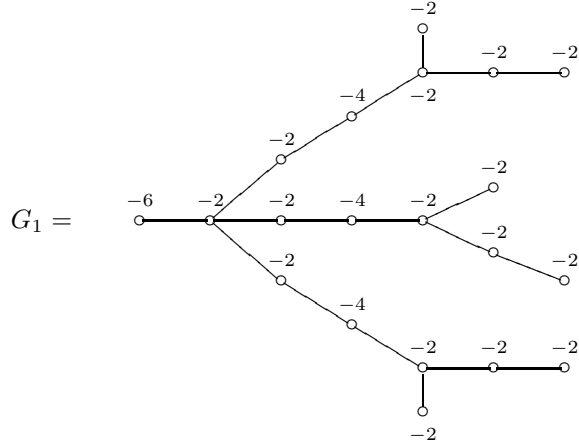
\xymatrix@R=6pt@C=24pt@M=0pt@W=0pt@H=0pt{
&&&&&&&\\
&&&&&\overtag{\Circ}{-2}{8pt}&&\\
&&&&&&&\\
&&&&&\undertag{\Circ}{-2}{4pt}\lineto[uu]\lineto[r] 
& \overtag{\Circ}{-2}{8pt}\lineto[r]&\overtag{\Circ}{-2}{8pt}
\\
&&&&&&&\\
&&&&\overtag{\Circ}{-4}{8pt}\lineto[uur]&&&\\
&&&&&&&\\
&&&\overtag{\Circ}{-2}{8pt}\lineto[uur]&&&&\\
&&&&&&
\overtag{\Circ}{-2}{8pt}&\\
 G_1=&
\overtag{\Circ}{-6}{8pt}\lineto[r] 
&\overtag{\Circ}{-2}{8pt}\lineto[uur] \lineto[r]
\lineto[ddr] &\overtag{\Circ}{-2}{8pt}\lineto[r]
&\overtag{\Circ}{-4}{8pt}\lineto[r]
& \overtag{\Circ}{-2}{8pt}
\lineto[ur]
\lineto[dr] &&\\
&&&&&& \overtag{\Circ}{-2}{8pt} \lineto[dr] &\\
&&&\overtag{\Circ}{-2}{8pt}\lineto[ddr]&&&&\overtag{\Circ}{-2}{8pt}\\
&&&&&&&\\
&&&&\overtag{\Circ}{-4}{8pt}\lineto[ddr] &&&\\
&&&&&&&\\
&&&&&\overtag{\Circ}{-2}{8pt}\lineto[dd]\lineto[r] & 
\overtag{\Circ}{-2}{8pt}\lineto[r] & \overtag{\Circ}{-2}{8pt}\\
&&&&&&&\\
&&&&&\undertag{\Circ}{-2}{4pt}&&\\
&&&&&&&\\
&&&&&&&\hbox to 0pt {~.\hss}}
$$  

If $\Gamma_{e_{N-2}}\neq\Gamma$ then one adds $G_1$ to the collection
of $\Delta_i$'s not used, and repeat the process by taking the to
plumbing diagrams of this collection which have arrows which triple
start with $N-2$. One continues this process until all the
$\Delta_i$'s have been used, and the final $G_N-1$ is then a plumbing
diagram for the universal abelian cover $\widetilde{M}$ of $M$. 

We will finish by giving a couple of other examples of the use of the
algorithm, but will leave the details of the calculation to the
readers.

\begin{ex}
Let $M$ be the manifold defined by the following plumbing diagram:
$$
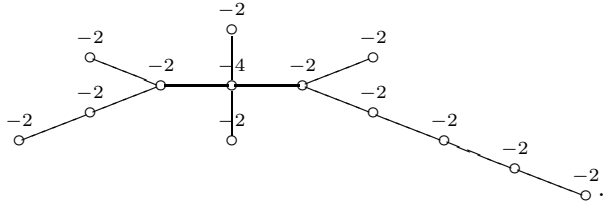
\xymatrix@R=6pt@C=24pt@M=0pt@W=0pt@H=0pt{
&&&\overtag{\Circ}{-2}{8pt}&&&&&\\
&\overtag{\Circ}{-2}{8pt}  &&
&&\overtag{\Circ}{-2}{8pt}&&&\\ 
&& \overtag{\Circ}{-2}{8pt}\lineto[ul]\lineto[dl]\lineto[r]
&\overtag{\Circ}{-4}{8pt}\lineto[r]\lineto[uu]\lineto[dd] &
\overtag{\Circ}{-2}{8pt}\lineto[ur]\lineto[dr]&&&&\\
&\overtag{\Circ}{-2}{8pt}\lineto[dl]&& & &
\overtag{\Circ}{-2}{8pt}\lineto[dr]&&&\\ 
\overtag{\Circ}{-2}{8pt}&&&\overtag{\Circ}{-2}{8pt}
&&&\overtag{\Circ}{-2}{8pt}\lineto[dr]&&\\ 
&&&&&&&\overtag{\Circ}{-2}{8pt}\lineto[dr]&\\
&&&&&&&&\overtag{\Circ}{-2}{8pt} \hbox to 0pt {~.\hss}}
$$
Its splice diagram is 
$$\splicediag{8}{30}{
  &&&&\Circ &&&\\
  &\Circ &&&&&&\Circ \\
  \Gamma=&&\overtag\Circ {v_1} {8pt}\lineto[ul]_(.5){2}
  \lineto[dl]^(.5){3}
  \lineto[rr]^(.25){44}^(.75){5}_(.5){e_1}&&\overtag\Circ {v_2\hspace{.4cm}}
  {8pt}\lineto[uu]_(.5){2} 
  \lineto[dd]^(.5){2}
  \lineto[rr]^(.25){7}^(.75){36}_(.5){e_2} &&
 \overtag\Circ{v_3}{8pt}
  \lineto[ur]^(.5){2}
  \lineto[dr]_(.5){5}&\\
  &\Circ &&&&&&\Circ \\
&&&&\Circ& &&\hbox to 0 pt{~.\hss} }$$
If we first cut along the edge called $e_1$ we get
$$\splicediag{8}{30}{
  &&&&&\Circ &&&\\
  &\Circ &&&&&&&\Circ \\
  \Gamma_{e_1}=&&\overtag\Circ {v_1} {8pt}\lineto[ul]_(.5){2}
  \lineto[dl]^(.5){3}
  \lineto[r]^(.5){22} & \overtag\Circ {(1,2)} {8pt}   & \overtag\Circ
  {(1,1)} {8pt}  &\overtag\Circ
  {v_2\hspace{.4cm}} 
  {8pt}\lineto[l]_(.5){5} \lineto[uu]_(.5){2} 
  \lineto[dd]^(.5){2}
  \lineto[rr]^(.25){7}^(.75){36}_(.5){e_2} &&
 \overtag\Circ{v_3}{8pt}
  \lineto[ur]^(.5){2}
  \lineto[dr]_(.5){5}&\\
  &\Circ &&&&&&&\Circ \\
&&&&&\Circ& &&\hbox to 0 pt{~.\hss} }$$
And cutting along $e_2$ the gives us
$$\splicediag{8}{30}{
  &&&&&\Circ &&&&\\
  &\Circ &&&&&&&&\Circ \\
  \Gamma_{e_2}=&&\overtag\Circ {v_1} {8pt}\lineto[ul]_(.5){2}
  \lineto[dl]^(.5){3}
  \lineto[r]^(.5){22} & \overtag\Circ {(1,2)} {8pt}   & \overtag\Circ
  {(1,1)} {8pt}  &\overtag\Circ
  {v_2\hspace{.4cm}} 
  {8pt}\lineto[l]_(.5){5} \lineto[uu]_(.5){2} 
  \lineto[dd]^(.5){2}
  \lineto[r]^(.5){7} &\overtag\Circ {(2,1)} {8pt} &\overtag\Circ
  {(2,2)} {8pt} & \overtag\Circ{v_3}{8pt}\lineto[l]_(.5){18}
  \lineto[ur]^(.5){2}
  \lineto[dr]_(.5){5}&\\
  &\Circ &&&&&&&&\Circ \\
&&&&&\Circ&& &&\\
&&\Gamma_1&&&\Gamma_2&&&\Gamma_3&\hbox to 0 pt{~.\hss} }$$
Next one determines the $3$ building blocks and get the following
plumbing diagrams
$$
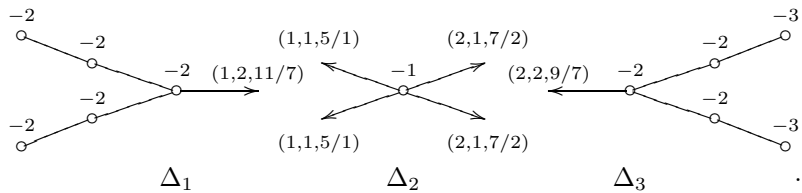
\xymatrix@R=6pt@C=24pt@M=0pt@W=0pt@H=0pt{
\overtag{\Circ}{-2}{8pt}&&&&&&&&&& \overtag{\Circ}{-3}{8pt}\\
&\overtag{\Circ}{-2}{8pt}\lineto[ul]  &&
&&&&&&\overtag{\Circ}{-2}{8pt}\lineto[ur]&\\ 
&& \overtag{\Circ}{-2}{8pt}\ar[r]^(1){(1,2,11/7)}
\lineto[ul]\lineto[dl] 
&&&\overtag{\Circ}{-1}{8pt}\ar[ul]_(1.25){(1,1,5/1)}
\ar[dl]^(1.25){(1,1,5/1)} \ar[ur]^(1.25){(2,1,7/2)}
\ar[dr]_(1.25){(2,1,7/2)} &&&
\overtag{\Circ}{-2}{8pt}\lineto[ur]\lineto[dr]\ar[l]_(1){(2,2,9/7)} &&\\
&\overtag{\Circ}{-2}{8pt}\lineto[dl]  &&
&&&&&&\overtag{\Circ}{-2}{8pt}\lineto[dr]&\\ 
\overtag{\Circ}{-2}{8pt}&&&&&&&&&& \overtag{\Circ}{-3}{8pt}\\
&&\Delta_1&&&\Delta_2&&&\Delta_3&&
 \hbox to 0pt {~.\hss}}
$$
One first glue the one copy of $\Delta_2$ to two copies of $\Delta_3$
and get after calculating the strings 
$$
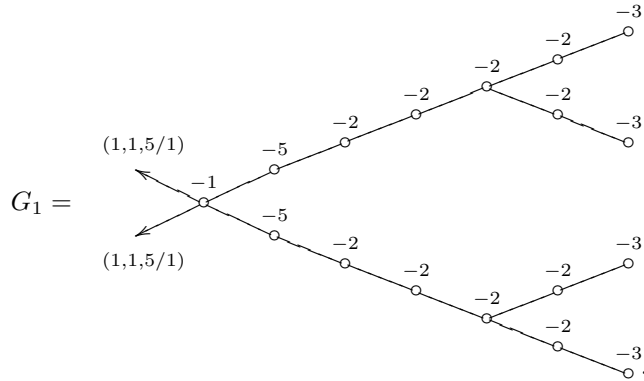
\xymatrix@R=6pt@C=24pt@M=0pt@W=0pt@H=0pt{
&&&&&&&& \overtag{\Circ}{-3}{8pt}\\
&&&&&&&\overtag{\Circ}{-2}{8pt}\lineto[ur]&\\
&&&&&&\overtag{\Circ}{-2}{8pt}\lineto[ur]\lineto[dr]&&\\ 
&&&&&\overtag{\Circ}{-2}{8pt}\lineto[ur] &&
\overtag{\Circ}{-2}{8pt}\lineto[dr]&\\ 
&&&&\overtag{\Circ}{-2}{8pt}\lineto[ur] &&&&\overtag{\Circ}{-3}{8pt}\\ 
&&&\overtag{\Circ}{-5}{8pt}\lineto[ur]&&&&&\\
G_1=&&\overtag{\Circ}{-1}{8pt}\ar[ul]_(1.25){(1,1,5/1)}
\ar[dl]^(1.25){(1,1,5/1)} \lineto[ur]
\lineto[dr] &&&&&&\\
&&&\overtag{\Circ}{-5}{8pt}\lineto[dr]&&&&&\\
&&&&\overtag{\Circ}{-2}{8pt}\lineto[dr] &&&&\overtag{\Circ}{-3}{8pt}\\
&&&&&\overtag{\Circ}{-2}{8pt}\lineto[dr] &&
\overtag{\Circ}{-2}{8pt}\lineto[ur]&\\
&&&&&&\overtag{\Circ}{-2}{8pt}\lineto[ur]\lineto[dr]&&\\ 
&&&&&&&\overtag{\Circ}{-2}{8pt}\lineto[dr]&\\
&&&&&&&& \overtag{\Circ}{-3}{8pt}
 \hbox to 0pt {~.\hss}}
$$
Then gluing two copies of $\Delta_1$ to $G_1$ and  calculate the
strings gives the following plumbing diagram for the universal abelian
cover
$$
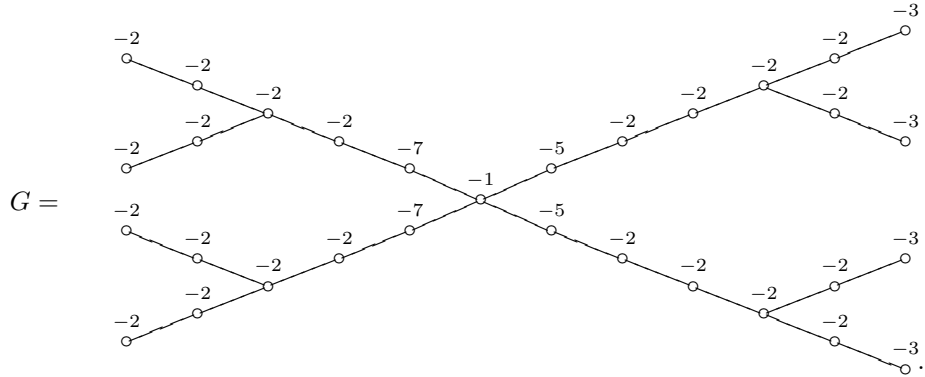
\xymatrix@R=6pt@C=24pt@M=0pt@W=0pt@H=0pt{
&&&&&&&&&&&& \overtag{\Circ}{-3}{8pt}\\
&\overtag{\Circ}{-2}{8pt}\lineto[dr]
&&&&&&&&&&\overtag{\Circ}{-2}{8pt}\lineto[ur]&\\ 
&&\overtag{\Circ}{-2}{8pt}\lineto[dr]
&&&&&&&&\overtag{\Circ}{-2}{8pt}\lineto[dr]\lineto[ur]&&\\  
&&&\overtag{\Circ}{-2}{8pt}\lineto[dr]
&&&&&&\overtag{\Circ}{-2}{8pt}\lineto[ur]  &&
\overtag{\Circ}{-2}{8pt}\lineto[dr]&\\
&&\overtag{\Circ}{-2}{8pt}\lineto[ur]&&\overtag{\Circ}{-2}{8pt}\lineto[dr]&
&&&\overtag{\Circ}{-2}{8pt}\lineto[ur] &&&&\overtag{\Circ}{-3}{8pt}\\
&\overtag{\Circ}{-2}{8pt}\lineto[ur]&&&& \overtag{\Circ}{-7}{8pt}\lineto[dr]
&&\overtag{\Circ}{-5}{8pt}\lineto[ur]&&&&\\
G=&&&&&&\overtag{\Circ}{-1}{8pt}\lineto[ur]\lineto[dr] &&&&&&\\
&\overtag{\Circ}{-2}{8pt}\lineto[dr]&&&& \overtag{\Circ}{-7}{8pt}\lineto[ur]
&&\overtag{\Circ}{-5}{8pt}\lineto[dr]&&&&\\
&&\overtag{\Circ}{-2}{8pt}\lineto[dr]&&\overtag{\Circ}{-2}{8pt}\lineto[ur]&
&&&\overtag{\Circ}{-2}{8pt}\lineto[dr] &&&&\overtag{\Circ}{-3}{8pt}\\
&&&\overtag{\Circ}{-2}{8pt}\lineto[ur]
&&&&&&\overtag{\Circ}{-2}{8pt}\lineto[dr]  &&
\overtag{\Circ}{-2}{8pt}\lineto[ur]&\\
&&\overtag{\Circ}{-2}{8pt}\lineto[ur]
&&&&&&&&\overtag{\Circ}{-2}{8pt}\lineto[ur]\lineto[dr]&&\\  
&\overtag{\Circ}{-2}{8pt}\lineto[ur]
&&&&&&&&&&\overtag{\Circ}{-2}{8pt}\lineto[dr]&\\ 
&&&&&&&&&&&& \overtag{\Circ}{-3}{8pt}
 \hbox to 0pt {~.\hss}}
$$

\end{ex}

\begin{ex}
Let $M$ be the graph manifold with the following plumbing diagram
$$
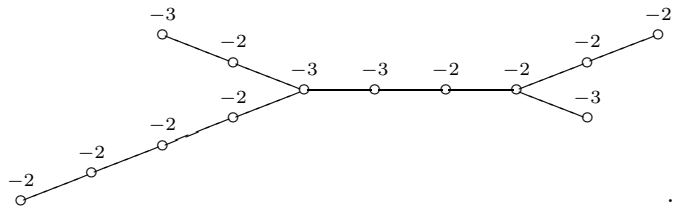
\xymatrix@R=6pt@C=24pt@M=0pt@W=0pt@H=0pt{
&&\overtag{\Circ}{-3}{8pt}&&&&&&&\overtag{\Circ}{-2}{8pt}\\
&&&\overtag{\Circ}{-2}{8pt}\lineto[ul]  &
&&&&\overtag{\Circ}{-2}{8pt}\lineto[ur]&\\ 
&&&& \overtag{\Circ}{-3}{8pt}\lineto[ul]\lineto[dl]\lineto[r]
&\overtag{\Circ}{-3}{8pt}\lineto[r]
&\overtag{\Circ}{-2}{8pt}\lineto[r] &
\overtag{\Circ}{-2}{8pt}\lineto[ur]\lineto[dr]&&\\
&&&\overtag{\Circ}{-2}{8pt}\lineto[dl]&&& & &
\overtag{\Circ}{-3}{8pt}&\\ 
&&\overtag{\Circ}{-2}{8pt}\lineto[dl] &&&&&&&\\ 
&\overtag{\Circ}{-2}{8pt}\lineto[dl]&&&&&&&&\\
\overtag{\Circ}{-2}{8pt}&&&&&&&&& \hbox to 0pt {~.\hss}}
$$
Its splice diagram the becomes
$$\splicediag{8}{30}{
  &&&&&\\
  &\Circ &&&&\Circ \\
  \Gamma=&&\overtag\Circ {v_0} {8pt}\lineto[ul]_(.5){5}
  \lineto[dl]^(.5){5}
  \lineto[rr]^(.25){27}^(.75){150}&& \overtag\Circ{v_1}{8pt}
  \lineto[ur]^(.5){3}
  \lineto[dr]_(.5){3}&\\
  &\Circ &&&&\Circ \\
&&&&&\hbox to 0 pt{~.\hss} }$$
Cutting the edge gives us the one node splice diagrams
$$\splicediag{8}{30}{
  &  \Circ &&&&&&\Circ \\
  \Gamma_1=&&\overtag\Circ {v_0} {8pt}\lineto[ul]_(.5){5}
  \lineto[dl]^(.5){5}
  \lineto[r]^(.25){9}& \overtag\Circ {(1,3)} {8pt},&\Gamma_2=&  \overtag\Circ
  {(1,35} {8pt}\lineto[r]^(.75){30}& \overtag\Circ{v_1}{8pt}
  \lineto[ur]^(.5){5}
  \lineto[dr]_(.5){5}&\\
  &\Circ &&&&&&\Circ\hbox to 0 pt{~,\hss} }$$
and the building blocks becomes
$$
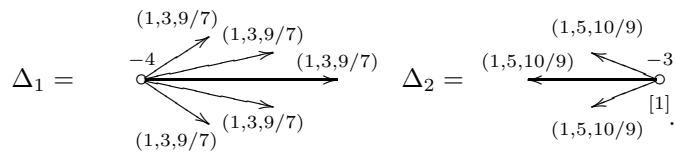
\xymatrix@R=6pt@C=24pt@M=0pt@W=0pt@H=0pt{
&&&&&&&& \\
&  &&
&&&&&\\ 
\Delta_1=& \overtag{\Circ}{-4}{8pt}\ar[uur]^(1){(1,3,9/7)}
\ar[urr]^(1){(1,3,9/7)} \ar[rrr]^(1){(1,3,9/7)}
\ar[drr]_(1){(1,3,9/7)} \ar[ddr]_(1){(1,3,9/7)} 
&&&&\Delta_2=&&& 
\undertag{\overtag{\Circ}{-3}{8pt}}{[1]}{4pt}\ar[ul]_(1.25){(1,5,10/9)}
\ar[dl]^(1.25){(1,5,10/9)} \ar[ll]_(1){(1,5,10/9)}\\
 &&
&&&&&\\ 
&&&&&&&& 
 \hbox to 0pt {~.\hss}}
$$
so to create the plumbing diagram $G$ of the universal abelian cover,
we glue $3$ copies of $\Delta_1$ to $5$ copies of $\Delta_2$ calculate
the string and get
$$
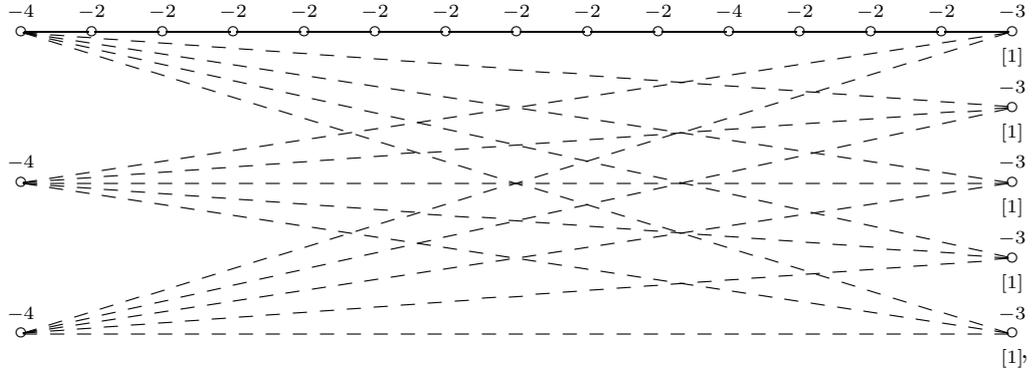
\xymatrix@R=6pt@C=24pt@M=0pt@W=0pt@H=0pt{
&&&&&&&&&&&&&&\\
\overtag{\Circ}{-4}{8pt}\lineto[r]\dashto[rrrrrrrrrrrrrrdddd]
\dashto[rrrrrrrrrrrrrrdddddddd] \dashto[rrrrrrrrrrrrrrdddddddddddd]
\dashto[rrrrrrrrrrrrrrdddddddddddddddd]
&\overtag{\Circ}{-2}{8pt}\lineto[r]&
\overtag{\Circ}{-2}{8pt}\lineto[r]&\overtag{\Circ}{-2}{8pt}\lineto[r]&
\overtag{\Circ}{-2}{8pt}\lineto[r]&\overtag{\Circ}{-2}{8pt}\lineto[r]&
\overtag{\Circ}{-2}{8pt}\lineto[r]&\overtag{\Circ}{-2}{8pt}\lineto[r]&
\overtag{\Circ}{-2}{8pt}\lineto[r]&\overtag{\Circ}{-2}{8pt}\lineto[r]&
\overtag{\Circ}{-4}{8pt}\lineto[r]&\overtag{\Circ}{-2}{8pt}\lineto[r]&
\overtag{\Circ}{-2}{8pt}\lineto[r]&\overtag{\Circ}{-2}{8pt}\lineto[r]&
\undertag{\overtag{\Circ}{-3}{8pt}}{[1]}{4pt}
\\
&&&&&&&&&&&&&&\\
&&&&&&&&&&&&&&\\
&&&&&&&&&&&&&&\\
&&&&&&&&&&&&&&\undertag{\overtag{\Circ}{-3}{8pt}}{[1]}{4pt} \\
&&&&&&&&&&&&&&\\
&&&&&&&&&&&&&&\\
&&&&&&&&&&&&&&\\
\overtag{\Circ}{-4}{8pt}\dashto[rrrrrrrrrrrrrruuuu]
\dashto[rrrrrrrrrrrrrruuuuuuuu] \dashto[rrrrrrrrrrrrrr]
\dashto[rrrrrrrrrrrrrrdddd] \dashto[rrrrrrrrrrrrrrdddddddd]
&&&&&&&&&&&&&&
\undertag{\overtag{\Circ}{-3}{8pt}}{[1]}{4pt}\\ 
&&&&&&&&&&&&&&\\
&&&&&&&&&&&&&&\\
&&&&&&&&&&&&&&\\
&&&&&&&&&&&&&&\undertag{\overtag{\Circ}{-3}{8pt}}{[1]}{4pt} \\
&&&&&&&&&&&&&&\\
&&&&&&&&&&&&&&\\
&&&&&&&&&&&&&&\\
\overtag{\Circ}{-4}{8pt}\dashto[rrrrrrrrrrrrrruuuu]
\dashto[rrrrrrrrrrrrrruuuuuuuu] \dashto[rrrrrrrrrrrrrruuuuuuuuuuuu]
\dashto[rrrrrrrrrrrrrruuuuuuuuuuuuuuuu] \dashto[rrrrrrrrrrrrrr]
&&&&&&&&&&&&&&
\undertag{\overtag{\Circ}{-3}{8pt}}{[1]}{4pt}\\ 
&&&&&&&&&&&&&& \hbox to 0pt {~,\hss}}
$$
where all the dashed lines represent strings identical to the string
at the top. Also remember that the graph is not a planar graph so any
intersection between the strings represented by the dashed lines do
not represent intersections in $G$, just crossings arising by a planar
projection of $G$ which is what we see here.
\end{ex}

\newpage

\bibliography{algorithm}

\end{document}